\documentclass[12pt]{article}
\usepackage[utf8]{inputenc}
\usepackage{amssymb,amsmath,amsfonts,amsthm,amscd,latexsym,indentfirst,verbatim}
\usepackage[T2A]{fontenc}
\usepackage{geometry}
\def\dbl{\lbrace\kern-3pt\lbrace}
\def\dbr{\rbrace\kern-3pt\rbrace}

\def\End{\textrm{End}}
\def\Hom{\textrm{Hom}}

\def\As{\textrm{As}}
\def\Rep{\textrm{Rep}}
\def\Span{\textrm{Span}}
\def\id{\textrm{id}}
\def\tr{\textrm{tr}}
\def\Spec{\mathrm{Spec}\,}

\begin{document}

\sloppy

\hfill{16W99, 17B38, 17B63 (MSC2020)}

\begin{center}
{\Large
Double Lie algebras of a nonzero weight}

\smallskip

Maxim Goncharov, Vsevolod Gubarev
\end{center}

\begin{abstract}
We introduce the notion of $\lambda$-double Lie algebra,
which coincides with usual double Lie algebra when $\lambda = 0$.
We show that every $\lambda$-double Lie algebra for $\lambda\neq0$ provides
the structure of modified double Poisson algebra on the
free associative algebra. In particular, it confirms the conjecture of S. Arthamonov (2017).
We prove that there are no simple finite-dimensional
$\lambda$-double Lie algebras.

\end{abstract}

\section{Introduction}

The notion of a~double Poisson algebra on a given associative algebra
was introduced by M. Van den Bergh in 2008~\cite{DoublePoisson}
as a noncommutative analog of Poisson algebra.
The goal behind this notion was to develop noncommutative Poisson geometry.
Let us briefly give the background of this object.

Given a finitely generated associative algebra~$A$ and $n\in\mathbb{N}$,
consider the representation space $\mathrm{Rep}_n(A) = \Hom(A,M_n(F))$,
where $F$ denotes the ground field.
We want to equip~$A$ with a~structure such that
$\mathrm{Rep}_n(A)$ is a~Poisson variety for every $n$.

For $a\in A$, we may consider the matrix-valued
function $a_{ij}$ on $\mathrm{Rep}_n(A)$,\,$1\leq i,j\leq n$. \linebreak
These functions generate the coordinate ring
$\mathcal{O}(\mathrm{Rep}_n(A))$ and
satisfy the relations
$(ab)_{ij} = \sum \limits_{k=1}^n a_{ik}b_{kj}$.
Thus, to define of a Poisson bracket $\{\cdot,\cdot\}$ on
$\mathrm{Rep}_n(A)$ one should know the value of 
$\{a_{ij},b_{kl}\}$ for all $a,b\in A$.
For this reason, M. Van den Bergh defined a~bilinear double bracket
$\dbl \cdot,\cdot\dbr\colon A\otimes A\to A\otimes A$
satisfying the analogs of anti-commutativity and Leibniz rule (valued in $A\otimes A$),
as well as Jacobi identity (valued in $A\otimes A\otimes A$).
An associative algebra equipped with such
a~double bracket is called {\bf double Poisson algebra}.

In~\cite{DoublePoisson}, it was shown that 
given a double Poisson algebra~$(A,\cdot,\dbl\cdot,\cdot\dbr)$, we have that
$\mathcal{O}(\mathrm{Rep}_n(A))$ is a Poisson algebra under the bracket
$$
\{a_{ij},b_{kl}\}
 = (\dbl a,b\dbr_{(1)})_{kj}(\dbl a,b\dbr_{(2)})_{il},
$$
where $\dbl a,b\dbr = \dbl a,b\dbr_{(1)}\otimes \dbl a,b\dbr_{(2)}$.

If one deals only with~$A$, then the product
$\{a,b\} = \dbl a,b\dbr_{(1)}\dbl a,b\dbr_{(2)}$
is a~derivation on its second argument and vanishes on commutators on its first one.
Moreover, $(A,\{\cdot,\cdot\})$ satisfies the Leibniz rule
$\{a,\{b,c\}\} = \{\{a,b\},c\} + \{b,\{a,c\}\}$
and $(A/[A,A],\{\cdot,\cdot\})$ is a Lie algebra.
In the terminology of  W. Crawley-Boevey~\cite{Crawley-Boevey},
every double Poisson algebra provides a $H_0$-Poisson structure.

N.~Iyudu, M.~Kontsevich, and Y.~Vlassopoulos~\cite{IKV}
showed that every double Poisson algebra appears as a particular part 
of a~pre-Calabi-Yau structure.
In~\cite{Kac15}, A.~De~Sole, V.~G.~Kac, and D.~Valeri introduced and studied
double Poisson vertex algebras.

In the middle of 2010s, S.~Arthamonov introduced a notion
of~{\bf modified double Poisson algebra}~\cite{Arthamonov0,Arthamonov}
with weakened versions of anti-commutativity and Jacobi identity.
This notion allowed S.~Arthamonov to study the Kontsevich system
and give more examples of $H_0$-Poisson structures arisen from double brackets.

The notion of {\bf double Lie algebra} naturally arose
directly from the definition of double Poisson algebra, 
it is a vector space~$V$ endowed with a double bracket
satisfying anti-commutativity and Jacobi identity mentioned above,
and we forget about associative product on~$V$.
As far as we know, T. Schedler was the first who defined such a notion clearly~\cite{Schedler}.
The importance of double Lie algebras is the following: every double Lie algebra
structure defined on a vector space~$V$ can be uniquely extended to
a~double Poisson algebra structure on the free associative algebra $\As\langle V\rangle$.
In~\cite{DoublePoissonFree}, A.~Odesskii, V.~Rubtsov, V.~Sokolov extended
linear and quadratic double Lie algebras defined on a $n$-dimensional vector
space to double Poisson algebras defined 
on the free $n$-generated associative algebra.

It is known that double Lie algebras on a finite-dimensional vector space~$V$ are in one-to-one correspondence with skew-symmetric Rota---Baxter operators of weight~0 on the matrix algebra~$M_n(F)$, where $n = \dim(V)$~\cite{DoubleLie,DoublePoissonFree,Schedler}.
Recall that a linear operator~$R$ defined on an algebra~$A$ is called
a~{\bf Rota---Baxter operator} (RB-operator, for short) of weight~$\lambda$, if
$$
R(x)R(y) = R( R(x)y + xR(y) + \lambda xy )
$$
for all $x,y\in A$.
This notion firstly appeared in the article~\cite{Tricomi} of F. Tricomi in 1951
and further was several times rediscovered~\cite{Baxter,BelaDrin82}, see
also the monograph of L.~Guo~\cite{GuoMonograph}.
To the moment, applications of Rota---Baxter operators
in symmetric polynomials, quantum field renormalization,
pre- and postalgebras, shuffle algebra, etc. were found~\cite{Aguiar00,Atkinson,FardThesis,GuoMonograph,Ogievetsky}.

Let us mention the bijection~\cite{Aguiar00,Unital,Schedler} between RB-operators of weight~0 on the matrix algebra $M_n(F)$ and solutions of the {\bf associative Yang---Baxter
equation} (AYBE)~\cite{Aguiar01,Polishchuk,Zhelyabin} on $M_n(F)$.
Recently this correspondence was established~\cite{AYBE-ext}
in the weighted case of both: RB-operators and AYBE~\cite{FardThesis}.

In~\cite{Double-0}, the correspondence between double Lie algebras and skew-symmetric RB-operators of weight~0 on the matrix algebra was extended for the infinite-dimensional case. 

In~\cite{DoubleLie}, M. Goncharov and P. Kolesnikov proved
that there are no simple finite-dimensional double Lie algebras.
The example of a~countable-dimensional simple double Lie algebra was found in~\cite{Double-0}.

We apply Rota---Baxter operators of nonzero weight on the matrix algebra
to define a weighted analog of double Lie algebras.
In this way, we show that a naive version of such a definition fails, see~\S3.
However, we define what is a $\lambda$-double Lie algebra for a~fixed~$\lambda\in F$.
Thus, in the finite-dimensional case
we extend the bijections
$$
\begin{matrix}
\mbox{double Lie} \\
\mbox{algebra}
\end{matrix}
\Longleftrightarrow
\begin{matrix}\mbox{skew-symmetric RB-operator} \\
\mbox{of weight\,0 on }M_n(F) \\
\end{matrix}
\Longleftrightarrow
\begin{matrix}\mbox{skew-symmetric solution} \\
\mbox{of AYBE on }M_n(F)
\end{matrix}
$$
for the weighted analogs of the objects as follows,
$$
\begin{matrix}
\lambda\mbox{-double Lie} \\
\mbox{algebra}
\end{matrix}
\Longleftrightarrow
\begin{matrix}\lambda\mbox{-skew-symmetric RB-operator} \\
\mbox{of weight\,} \lambda\mbox{ on }M_n(F) \\
\end{matrix}
\Longleftrightarrow
\begin{matrix}(-\lambda)\mbox{-skew-symmetric solution} \\
\mbox{of AYBE}(-\lambda) \mbox{ on }M_n(F)
\end{matrix}
$$

The correspondence between $\lambda$-double Lie algebras
and RB-operators of weight~$\lambda$ is helpful for 
constructing examples of $\lambda$-double Lie algebras.
As in the case $\lambda = 0$, we prove that there are no simple
finite-dimensional $\lambda$-double Lie algebras.
Recall~\cite{Double-0} that a double Lie algebra $V$ is said to be simple if 
$\dbl V,V\dbr \neq (0)$ and there are no nonzero proper subspaces $I$ in $V$
such that $\dbl V,I\dbr + \dbl I,V\dbr \subseteq I\otimes V + V\otimes I$.

On the other hand, we find a~pair of interesting infinite-dimensional
$\lambda$-double Lie algebras, one of them is $F[t]$ equipped
with the double $\lambda$-skew-symmetric bracket
$$
\dbl t^n,t^m\dbr = \frac{t^m\otimes t^{n+1}-t^n\otimes t^{m+1}}{t\otimes 1-1\otimes t}.
$$
We show that this double Lie algebra~$M$ has the only one nonzero proper ideal
which occurs to be isomorphic to~$M$.

Finally, we prove that every $\lambda$-double Lie algebra structure on
a vector space~$V$ generates a~unique modified double Poisson algebra structure
on $\As\langle V\rangle$.
This general result confirms Conjecture~21 of S. Arthamonov (2017)~\cite{Arthamonov} about the double bracket~$\dbl\cdot,\cdot\dbr^{II}$
defined on the three-dimensional vector space $V = \Span\{a_1,a_2,a_3\}$ as follows,
\begin{gather*}
\dbl a_1,a_2\dbr^{II}  = -a_1\otimes a_2,\quad
\dbl a_2,a_1\dbr^{II}  =  a_1\otimes a_2,\quad
\dbl a_2,a_3\dbr^{II}  =  a_3\otimes a_2,\\
\dbl a_3,a_1\dbr^{II}  =  a_1\otimes a_3 - a_3\otimes a_1,\quad
\dbl a_3,a_2\dbr^{II}  = -a_3\otimes a_2.
\end{gather*}
Conjecture says that the double bracket $\dbl\cdot,\cdot\dbr^{II}$ can be extended to a modified double Poisson algebra structure
on $\As\langle a_1,a_2,a_3\rangle$.

We want to emphasize the following interesting parallelism.
In~1982, A.A. Belavin and V.G.  Drinfel'd proved~\cite{BelaDrin82} that
given a skew-symmetric solution $r = \sum\limits a_i\otimes b_i$ of the classical Yang---Baxter equation (CYBE) on a semisimple finite-dimensional Lie algebra~$L$,
we get a Rota---Baxter operator~$R$ of weight~0 on $L$ defined by the formula
$R(x) = \sum \langle a_i,x\rangle b_i$.
Here $\langle \cdot,\cdot\rangle$ denotes the Killing form on $L$.
In~2017, M. Goncharov~\cite{Goncharov2} proved that a solution of modified (i.\,e., with weakened skew-symmetricity) CYBE on a simple finite-dimensional Lie algebra $L$ gives by the same formula an RB-operator of a~nonzero weight on~$L$.

In the case of double algebras, skew-symmetric RB-operators of weight~0 on the matrix algebra
produce double Poisson algebras. On the other hand, $\lambda$-skew-symmetric RB-operators of nonzero weight~$\lambda$ on the matrix algebra give rise to modified double Poisson algebras,
structures with weaker anti-commutativity and Jacobi identity.

Let us give a short outline of the work.
In~\S2, we give required preliminaries on Rota---Baxter operators,
associative Yang---Baxter equation, and double Lie algebras (including 
infinite-dimensional case).

In~\S3, we show that a naive version of $\lambda$-double Lie algebra
fails because of the properties of RB-operators of nonzero weight
on the matrix algebra.

In~\S4, we give the main definition of $\lambda$-double Lie algebra.
We provide both finite-dimensional and infinite-dimensional examples
of $\lambda$-double Lie algebras including simple infinite-dimensional ones.
In the finite-dimensional case, we prove that there are no simple $\lambda$-double Lie algebras.

In~\S5, we prove that every $\lambda$-double Lie algebra on a vector space~$V$
can be uniquely extended to a modified double Poisson algebra on the free associative algebra
$\As\langle V\rangle$.

\section{Preliminaries}

\subsection{Rota---Baxter operators}

A~linear operator~$R$ defined on an (not necessarily associative) algebra~$A$ is called a~Rota---Baxter operator (RB-operator, for short) of weight~$\lambda$, if
$$
R(x)R(y) = R( R(x)y + xR(y) + \lambda xy )
$$
holds for all $x,y\in A$.

It is well-known that given an RB-operator of weight~$\lambda$ on an algebra $A$, we have that $\widetilde{R} = -R-\lambda\id$ is again an RB-operator of weight~$\lambda$ on $A$.

{\bf Proposition 1}~\cite{Unital}.
Let $A$ be an algebra, let $R$ be an RB-operator of weight~$\lambda$ on~$A$,
and let $\psi$ be either an automorphism or an antiautomorphism of $A$.
Then the operator $R^{(\psi)} = \psi^{-1}R\psi$ is an RB-operator of weight~$\lambda$ on~$A$.

In~\cite{Spectrum}, it was proved the following general property of RB-operators on unital algebras.

{\bf Theorem 1}.
Let $A$ be a finite-dimensional unital algebra over a field~$F$.
Given a~Rota---Baxter operator $R$ of weight~$\lambda$ on $A$,
we have $\Spec(R)\subset\{0,-\lambda\}$.

{\bf Corollary 1}.
Given a~Rota---Baxter operator~$R$ of nonzero weight~$\lambda$ on~$M_n(F)$, we have the decomposition
$A = \ker(R^{N})\oplus \ker(R+\lambda\id)^{N}$ (as subalgebras), where $N=n^2$.

Given an algebra~$A$ and an ideal~$J$ of~$A$,
a linear map $R\colon J\to A$ is called a Rota---Baxter operator
of weight~$\lambda\in F$ from~$J$ to~$A$ if
$$
R(a)R(b) = R(R(a)b + aR(b) + \lambda ab)
$$
for all $a,b\in J$. For $J = A$, we obtain the usual notion of an RB-operator
of weight~$\lambda$ on~$A$.

\subsection{Associative Yang---Baxter equation}

Let $A$ be an associative algebra, $r = \sum a_i\otimes b_i\in A\otimes A$.
The tensor $r$ is a solution of the associative Yang---Baxter equation
(AYBE, \cite{Aguiar01,Polishchuk,Zhelyabin}) if
\begin{equation}\label{AYBE}
r_{13}r_{12}-r_{12}r_{23}+r_{23}r_{13} = 0,
\end{equation}
where 
$$
r_{12} = \sum a_i\otimes b_i\otimes 1,\quad
r_{13} = \sum a_i\otimes 1\otimes b_i,\quad
r_{23} = \sum 1\otimes a_i\otimes b_i
$$
are elements from $A^{\otimes3}$.

The switch map $\tau\colon A\otimes A\to A\otimes A$
acts in the following way: $\tau(a\otimes b) = b\otimes a$.
The solution $r$ of AYBE is called skew-symmetric if $r + \tau(r) = 0$.

{\bf Proposition 2}~\cite{Aguiar00}.
Let $r = \sum a_i\otimes b_i$ be a solution of AYBE
on an associative algebra~$A$.
A linear map $P_r\colon A\to A$ defined as
\begin{equation}\label{AYBE2RB}
P_r(x) = \sum a_i x b_i
\end{equation}
is an RB-operator of weight zero on $A$.

Later, in 2006 K.~Ebrahimi-Fard defined in his Thesis~\cite[p.~113]{FardThesis}
the associative Yang---Baxter equation of weight~$\lambda$.
Given an associative algebra~$A$ and a tensor $r\in A\otimes A$,
we say that $r$~is a~solution of the associative Yang---Baxter equation
of weight~$\lambda$ if
\begin{equation}\label{wAYBE}
r_{13}r_{12}-r_{12}r_{23}+r_{23}r_{13} = \lambda r_{13}.
\end{equation}

{\bf Proposition 3}~\cite{FardThesis,AYBE-ext}.
Let $r = \sum a_i\otimes b_i$ be a solution of AYBE of weight~$\lambda$ on an associative algebra~$A$.
A linear map $P_r\colon A\to A$ defined by~\eqref{AYBE2RB} is an RB-operator of weight~$-\lambda$ on~$A$.

{\bf Theorem 2}~\cite{Unital,AYBE-ext}.
The map $r\to P_r$ is a bijection between the set of the solutions of AYBE of weight~$\lambda$ on $M_n(F)$ and the set of RB-operators of weight~$-\lambda$ on $M_n(F)$.

\subsection{Double Poisson and double Lie algebras}

Let $V$ be a linear space over $F$. Given $u\in V^{\otimes n}$ and $\sigma \in S_n$, $u^\sigma$ denotes the permutation of tensor factors.
By a double bracket on $V$ we call a linear map from $V\otimes V$ to $V\otimes V$.
Given an associative algebra~$A$, we consider the outer bimodule action of $A$ on $A\otimes A$:
$b(a\otimes a') c = (ba)\otimes (a'c)$.

{\bf Definition 1}~\cite{DoublePoisson}.
A double Poisson algebra is an associative algebra $A$ equipped with a~double bracket
satisfying the following identities for all $a,b,c\in A$
\begin{gather}
\dbl a,b\dbr =- \dbl b,a\dbr ^{(12)}, \label{antiCom} \\
\dbl a, \dbl b,c\dbr \dbr _L -\dbl b, \dbl a,c\dbr \dbr _R
 = \dbl \dbl a,b\dbr,c\dbr _L,  \label{Jacobi} \\
\dbl a,bc\dbr = \dbl a,b\dbr c + b\dbl a,c\dbr, \label{Leibniz}
\end{gather}
where
$\dbl a, b\otimes c \dbr _L = \dbl a,b \dbr \otimes c$,
$\dbl a, b\otimes c\dbr _R = (b\otimes \dbl a,c\dbr )$, and
$\dbl a\otimes b, c\dbr _L = (\dbl a,c\dbr \otimes b)^{(23)}$.
Anti-commutativity~\eqref{antiCom} and Leibniz rule~\eqref{Leibniz}
imply~\cite{DoublePoisson} the following equality
\begin{equation} \label{LeibnizTwo}
\dbl ab,c\dbr = a*\dbl b,c\dbr + \dbl a,c\dbr * b
\end{equation}
for the inner bimodule action of $A$ on $A\otimes A$:
$b*(a\otimes a')*c = (ac)\otimes (ba')$.

{\bf Definition 2}~\cite{DoublePoissonFree,Schedler,Kac15}.
A double Lie algebra is a linear space $V$ equipped with a~double bracket
satisfying the identities~\eqref{antiCom} and~\eqref{Jacobi}.

An {\it ideal} of a double Lie algebra $V$ is a subspace $I\subseteq V$ such that
$\dbl V,I\dbr + \dbl I,V\dbr \subseteq I\otimes V + V\otimes I$.
Given an ideal $I$ of a double Lie algebra~$V$, we have
a~natural structure of a double Lie algebra on the space $V/I$, i.\,e.,
$\dbl x + I, y + I \dbr = \dbl x,y \dbr + I\otimes V + V\otimes I$.

Let~$L$ and $L'$ be double Lie algebras and let $\varphi\colon L\to L'$
be a~linear map. Then~$\varphi$ is called as {\it homomorphism} from~$L$ to~$L'$ if
$$
(\varphi\otimes\varphi)(\dbl a,b\dbr) = \dbl \varphi(a),\varphi(b)\dbr
$$
holds for all $a,b\in L$.
Note that the kernel of any homomorphism from~$L$ is an ideal of~$L$.
A bijective homomorphism from~$L$ to $L'$ is called an isomorphism of double Lie algebras.

A double Lie algebra $V$ is said to be {\it simple} if $\dbl V,V\dbr \neq (0)$
and there are no nonzero proper ideals in $V$.

Suppose $V$ is a finite-dimensional space.
In~\cite{DoubleLie}, it was shown that every double bracket 
$\dbl \cdot,\cdot\dbr$ on~$V$ is determined by a linear operator $R\colon\End (V)\to \End(V)$, precisely,
\begin{equation}\label{eq:Bracket_via_RB}
 \dbl a,b\dbr  = \sum\limits_{i=1}^N e_i(a)\otimes R(e_i^*)(b)
 , \quad a,b\in V,
\end{equation}
where $e_1,\dots, e_N$ is a linear basis of $\End(V)$, $e_1^*,\dots, e_N^*$
is the corresponding dual basis relative to the trace form.

Let us explain how to relate an operator $R$ on $\End(V)$ with a~bracket $\dbl\cdot,\cdot\dbr$ explicitly. 
Let $f_1,\ldots,f_n$ be a basis of $V$. By $e_{kl}$, $k,l\in\{1,\ldots,n\}$,
we mean the standard basis of~$\End(V)$, which acts on~$V$ by the formula
$e_{ij}(f_k) = \delta_{jk} f_i$. 
Let us rewrite
$$
\dbl f_k,f_l\dbr = \sum\limits_{m=1}^n f_m\otimes v_m,\quad k,l=1,\ldots,n,
$$
for some $v_m\in V$.
Define an operator $R\in\End(\End(V))$ as follows, $R(e_{km})(f_l) = v_m$,
where $k,l,m=1,\ldots,n$. Then 
$$
\dbl f_k,f_l\dbr 
 = \sum\limits_{m=1}^n f_m\otimes v_m
 = \sum\limits_{m,p=1}^n e_{mp}(f_k)\otimes R(e_{mp}^*)(f_l).
$$
By linearity, we obtain the formula~\eqref{eq:Bracket_via_RB} for all $a,b\in V$. 

Note that the identity operator corresponds to the switch map $\tau\colon f_k\otimes f_l\to f_l\otimes f_k$. 

A linear operator $P$ on~$\End(V)$ is called skew-symmetric if $P = -P^*$,
where $P^*$ is the conjugate operator on $\End(V)$ relative to the trace form.

{\bf Theorem~3}~\cite{DoubleLie}.
Let $V$~be a finite-dimensional vector space with a double bracket $\dbl \cdot,\cdot \dbr$
determined by an operator $R\colon\End(V)\to \End(V)$ by~\eqref{eq:Bracket_via_RB}.
Then $V$ is a double Lie algebra if and only if $R$ is a~skew-symmetric RB-operator of weight~0 on $\End(V)$.

{\bf Remark 1}. 
Theorem~3 was stated in~\cite{Schedler} in terms of skew-symmetric solutions of the associative Yang---Baxter equation (AYBE). Since there exists the one-to-one correspondence between
solutions of AYBE and Rota---Baxter operators of weight~0 on the matrix algebra~\cite{Unital},
Theorem~3 follows from~\cite{Schedler}. Actually, Theorem~3 was mentioned also in~\cite{DoublePoissonFree}.

\subsection{Infinite-dimensional double Lie algebras}

Consider a countable-dimensional double Lie algebra $\langle V,\dbl \cdot,\cdot\dbr\rangle$.
We fix a~linear basis $u_i$, $i\in\mathbb{N}$, of $V$.
Define $e_{ij}\in\End(V)$ by the formula $e_{ij}u_k = \delta_{jk}u_i$.
Let $\varphi\in\End(V)$, then we may write $\varphi = \sum\limits_{ij}a_{ij}e_{ij}$.
We identify $\varphi$ with an infinite matrix $[\varphi] = (a_{ij})_{i,j\geq0}$.
Since $\varphi\in\End(V)$ is well-defined,
there is only finite number of nonzero elements in every column of the matrix $[\varphi]$.

Define the subalgebra~$\End_f(V)$ of $\End(V)$ as follows,
$$
\End_f(V) = \{\varphi\in\End(V)\mid
 \mbox{ for every }i,\ [\varphi]_{ij} = 0
 \mbox{ for almost all }j\}.
$$

Denote by $I$ a linear span of matrix unities $e_{ij}$, it is an ideal in $\End_f(V)$.

Let $\varphi\in\End_f(V) = \sum\limits_{i,j}a_{ij}e_{ij}$.
We define the symmetric non-degenerate bilinear trace form
$\langle \cdot,\cdot \rangle$ on $I\times \End_f(V)\cup \End_f(V)\times I$ as follows,
$$
\langle e_{kl},\varphi \rangle
 = \langle \varphi,e_{kl}\rangle
 = \tr(e_{kl}\varphi)
 = a_{lk}.
$$
Moreover, the form is associative, i.\,e., $\langle a,bc\rangle = \langle ab,c\rangle$,
where at least one of $a,b,c$  lies in~$I$ and others are from $\End_f(V)$.

Given a double bracket algebra~$\dbl \cdot,\cdot\dbr$ on a space~$V$, 
we may define a linear operator $R\colon I\to \End (V)$ by the formula
\begin{equation}\label{Bracket_via_RBInf}
\dbl a,b\dbr  = \sum\limits_{i,j\geq0} e_{ij}(a)\otimes R(e_{ji})(b), \quad a,b\in V.
\end{equation}
Conversely, given an operator $R\colon I\to \End (V)$, 
one can define a double bracket on~$V$ by the formula~\eqref{Bracket_via_RBInf}.

We may define a conjugate operator $R^*\colon I\to \End (V)$ as follows,
\begin{equation}\label{Conjugate}
\dbl b,a\dbr^{(12)}  = \sum\limits_{i,j\geq0} e_{ij}(a)\otimes R^*(e_{ji})(b), 
\quad a,b\in V.
\end{equation}

In~\cite{Double-0}, the following equality was shown,
\begin{equation}\label{InfConj}
\langle R(x),y\rangle = \langle x,R^*(y)\rangle,\quad x,y\in I.
\end{equation}

{\bf Theorem~4}~\cite{Double-0}.
Let $V$~be a countable-dimensional vector space with a fixed linear basis
$u_i$ and with a~double bracket $\dbl \cdot,\cdot \dbr$
determined by a linear map $R\in \End_f(V)$ as in~\eqref{Bracket_via_RBInf}.
Then $V$ is a double Lie algebra if and only if $R$ is a~skew-symmetric RB-operator of weight~0 from~$I$ to $\End(V)$.

{\bf Theorem 5}~\cite{Double-0}.
The double Lie algebra $L_2$ defined on $F[t]$ by the formula
$$
\dbl t^n,t^m\dbr = \frac{t^n\otimes t^m-t^m\otimes t^n}{t\otimes 1-1\otimes t}
$$
is simple.

\section{Naive version of $\lambda$-double Lie algebras}

Let us return to the finite-dimensional case.
Suppose that we have a double bracket $\dbl\cdot,\cdot\dbr$ on a linear space $V$.
We want to apply the connection~\eqref{eq:Bracket_via_RB} of the double bracket
with a linear operator $R$ on $\End (V)$. By Theorem~3, the double bracket is Lie if and only if $R$ is a skew-symmetric RB-operator of weight~0 on $\End (V)$.

What does happen if $R$ is an RB-operator of nonzero weight~$\lambda$ on $\End (V)$? We have the Jacobi identity~\eqref{Jacobi} if and only if (see the proofs of Theorems~3 and~4)
$$
R(\theta_R(y)x) = 0,\quad x,y\in \End(V),
$$
here $\theta_R = R^* + R + \lambda\id$.
Thus, to avoid degenerate~$R$, we need 
$$\theta_R =R+R^*+\lambda \id= 0$$
which is some kind of skew-symmetricity in the nonzero weight case.
Hence, we have~\eqref{Jacobi} and the following identity,
\begin{equation}\label{bad-anticommutativity}
\dbl a,b\dbr
 + \dbl b,a\dbr^{(12)} = - \lambda b\otimes a,
\end{equation}
which is an analogue of anticommutativity.

Therefore, we want to define $\lambda$-analogue of double Lie algebra as a vector space~$V$ with a~double bracket satisfying~\eqref{Jacobi} and~\eqref{bad-anticommutativity}. 

Let us prove that such a notion is not natural. More precisely, we show that such finite-dimensional objects do not exist for $\lambda \neq 0$, and we do not expect that they exist in the infinite-dimensional case. 

By a variety of algebras we mean a class of all algebras satisfying the prescribed set of identities. 
For instance, the variety of associative algebras is defined by the identity $(xy)z = x(yz)$. 
By the famous Birkhoff theorem, a class~$\mathcal{M}$ of algebras forms a variety if and only if $\mathcal{M}$ is closed under taking of homomorphic images, subalgebras and direct products.

Let $\mathcal{M}$ be a variety of algebras, $A$ be an  
algebra from $\mathcal{M}$ over a field $F$, $R\colon A\mapsto A$ be a linear map, and $\lambda\in F$. Consider a~direct sum of vector spaces
$D_R(A)=A\oplus \bar{A}$, where $\bar{A}$ is an isomorphic copy of $A$. Define a product on $D_R(A)$ as follows
\begin{multline*}
(a+\bar{b})*(x+\bar{y}) \\
 = ax+R(ay)-aR(y)+R(bx)-R(b)x+\overline{ay+bx-R(b)y-bR(y)-\lambda by}.
\end{multline*}

{\bf Definition 3}. 
A bilinear non-degenerate symmetric form $\omega:A\times A\mapsto F$ on an algebra $A$ is called invariant if for all $a,b,c\in A$:
$$
\omega(ab,c)=\omega(a,bc).
$$
In this case, the pair $(A,\omega)$ is called a quadratic algebra.

The most important exampels of quadratic algebras are:

1. A semisimple finite-dimensional Lie algebra over a field of characteristic zero with the Killing form.

2. A matrix algebra $M_n(F)$ with the form $\omega:M_n(F)\times M_n(F)\mapsto F$ defined as
$$
\omega(a,b)=\tr{(ab)}.
$$

{\bf Remark 2}. 
If $(A,\omega)$ is a quadratic algebra, then $D_R(A)$ is isomorphic to the Drinfeld double of the bialgebra $(A,\delta_r)$ (here we mean bialgebra in general sense, that is, an algebra with a comultiplication), where $r=\sum a_i\otimes b_i\in A\otimes A$ corresponds to the map $R$ by the natural isomorphism of $\End(A)$ and $A\otimes A$, and 
$\delta_r(a) = \sum (a_ia\otimes b_i-a_i\otimes ab_i)$ (see \cite{Goncharov2} for details).

{\bf Proposition 4}. 
Let $A$ be an algebra from a variety $\mathcal{M}$.
If $R$ is a Rota---Baxter operator of weight $\lambda$ on~$A$, 
then $D_R(A)$ is an algebra from $\mathcal{M}$.

{\sc Proof}.
Set 
\begin{equation} \label{i-map}
i(a)=\bar{a}+(\lambda a+R(a)). 
\end{equation}
Let us prove that
$I(A)=\{i(a)\mid a\in A\}$ is an ideal of $D_R(A)$. 
Clearly, $I(A)$ is a subspace of $D_R(A)$  and $\dim(I(A)) = \dim A$. Let $a,b\in A$.
Then
\begin{equation}\label{u1}
i(a)*b
 = (\bar{a}+\lambda a+R(a))*b
 = \overline{ab}+R(ab)-R(a)b+\lambda ab+R(a)b
 = i(ab),
\end{equation}

\vspace{-0.9cm}
\begin{multline}\label{u2}
i(a)*\bar{b}
 = (\bar{a}+\lambda a+R(a))*\bar{b} \\
 = \overline{-\lambda ab-aR(b)-R(a)b+\lambda ab+R(a)b} \\
   +R(\lambda ab+R(a)b)-\lambda aR(b)-R(a)R(b) \\
 = \overline{-aR(b)}-\lambda aR(b)-R(aR(b))=i(-aR(b)).
\end{multline}
Hence, $I(A)$ is a left ideal of $D_R(A)$. 
Similarly, $I(A)$ is a right ideal of $D_R(A)$ and therefore,
$I(A)$  is an ideal of $D_R(A)$.

Moreover,
\begin{equation}\label{u3}
i(a)*i(b)
 = i(a)*(\bar{b}+\lambda b+R(b))
 = -i(aR(b))+i(\lambda ab+aR(b))
 = \lambda i( ab).
\end{equation}
If $\lambda\neq 0$, then a map
$i_{\lambda}\colon a\mapsto \frac{i(a)}{\lambda}$ is an isomorphism of algebras $A$ and $I(A)$. If $\lambda=0$, then $I(A)^2=0$. 
In both cases, the space $D_R(A)$ is equal to the direct sum of subalgebras~$A$ and $I(A)$. 
 
Now we have two situations: if $\lambda\neq 0$, then $D_R(A)$ is isomorphic to
$A\otimes D$, where $D$ is two-dimensional $F$-algebra with a basis  $\{1,d\}$ subject to the relation  $d^2=d$. 

If $\lambda=0$, then $D_R(A)$ is isomorphic to
$A\otimes N$, where $N$ is two-dimensional algebra with a basis $\{1,n\}$ subject to the relation $n^2=0$. Hence, $D_R(A)\in \mathcal{M}$. 
\hfill $\square$

{\bf Remark 3}. 
Let $\lambda\neq 0$, define the map $j\colon a\mapsto -\bar{a}-R(a)$ and  
\begin{equation} \label{J-map}
J(A)=\{j(a) \mid a\in A\}.
\end{equation}
Similar arguments as in Proposition 4 show that $J(A)$ is an ideal in $D_R(A)$ and the map $j_\lambda\colon a\mapsto \frac{j(a)}{\lambda}$ is an isomorphism of $A$ and $J(A)$. Note that  for all $a\in A$,  $j(a)+i(a)=\lambda a$. That is, $D_R(A)=I(A)\oplus J(A)$.

{\bf Remark 4}. 
A construction very close to $D_R(A)$ was suggested by K.~Uchino~\cite{Uchino}
when $A$ is associative and $\lambda = 0$.

Suppose in addition that a non-degenerate invariant symmetric bilinear form $\omega$ is defined on $A$. 
Then $\omega$ induces a form $Q$ on $D_R(A)$: for $a,b,c,d\in A$ put
\begin{equation}\label{f1}
Q(a+\bar{b},c+\bar{d})=\omega(a,d)+\omega(b,c).
\end{equation}

It is easy to see that $Q$ is a non-degenerate symmetric bilinear form on $D_R(A)$.

{\bf Proposition 5}. 
Let $(A,\omega)$ be a quadratic algebra and let $Q$ be the form defined on $D_R(A)$ by~\eqref{f1}. Then the form $Q$ is invariant if and only if for all $a,b\in A$:
\begin{equation}\label{p5.1}
    R(ab)+R^*(ab)+\lambda ab=0.
 \end{equation}

{\sc Proof}.
First of all, let us note that the following conditions are equivalent due to the non-degeneracy of the form on $A$:

1. For all $a,b\in A$: $R(ab)+R^*(ab)+\lambda ab=0$,

2. For all $a,b\in A$: $aR(b)+aR^*(b)+\lambda ab=0$,

3. For all $a,b\in A$: $R(a)b+R^*(a)b+\lambda ab=0$.

Indeed, let $a,b,c\in A$. Then
$$
\omega(R(ab)+R^*(ab)+\lambda ab,c)=\omega(a,bR^*(c)+bR(c)+\lambda bc).
$$
This shows the equivalence of 1 and 2. Similarly, 1 is equivalent to 3.

Let $a,b,c\in A$. From the definition of $Q$ we have that  $Q(A,A)=Q(\bar{A},\bar{A})=0$. 
Since $\omega$ is invariant, we have
$$
Q(a*b,\bar{c})=\omega(ab,c)=\omega(a,bc)=Q(a,\overline{bc})=Q(a,b*\bar{c}).
$$

Similarly, we can prove the following identities:
$$
Q(a*\bar{b},c)=Q(a,\bar{b}*c),\quad Q(\bar{a}*b,c)=Q(\bar{a},b*c).
$$

Further,
$$
Q(\bar{a}*\bar{b},c)=-\omega(R(a)b+aR(b)+\lambda ab,c)=-\omega(a,R^*(bc)+R(b)c+\lambda bc).
$$

On the other hand,
$$
Q(\bar{a},\bar{b}*c)=Q(\bar{a},\overline{bc}+R(bc)-R(b)c)=\omega(a,R(bc)-R(b)c).
$$

Therefore,
$$Q(\bar{a}*\bar{b},c)-Q(\bar{a},\bar{b}*c)=-\omega(a,R^*(bc)+R(bc)+\lambda bc).
$$
Since $\omega$ is non-degenerate, $Q(\bar{a}*\bar{b},c)-Q(\bar{a},\bar{b}*c)=0$  if and only if 
$$
R^*(bc)+R(bc)+\lambda bc=0
$$
for all $b,c\in A$. 

Similar arguments show that an identity $Q(a*\bar{b},\bar{c})=Q(a,\bar{b}*\bar{c})$ is also equivalent to the condition \eqref{p5.1}.

Consider equality $Q(\bar{a}*b,\bar{c})=Q(\bar{a},b*\bar{c})$. We have:
$$
Q(\bar{a}*b,\bar{c})=\omega(R(ab)-R(a)b,c)=\omega(a,bR^*(c)-R^*(bc)).
$$

Similarly,
$$
Q(\bar{a},b*\bar{c})=\omega(a,R(bc)-bR(c)).
$$

Therefore, $Q(\bar{a}*b,\bar{c})-Q(\bar{a},b*\bar{c})=0$ if and only if
\begin{equation}\label{p5.2}
R(bc)+R^*(bc)-bR(c)-bR^*(c)=0.
\end{equation}

It is only left to note that \eqref{p5.2} follows by \eqref{p5.1} and the observation from the beginning of the proof.
\hfill$\square$

{\bf Theorem 6}. 
Let $A=M_n(F)$ and let $\omega(x,y)=\tr(xy)$ be the trace form on $A$. 
There are no Rota---Baxter operators of a~nonzero weight~$\lambda$ on~$A$ satisfying the equality $R+R^*+\lambda \id=0$.

{\sc Proof}. 
Assume the converse. It is enough to consider the case when $\lambda=1$. By Proposition 5, $Q$ is a non-degenerate invariant bilinear form on $D_R(A)$. Let $E\in A$ be the identity matrix. Since $\tr(E)\neq 0$, then $Q(E,\bar{E})\neq 0$. 

Let $a,b\in A$. Then
\begin{equation}\label{th6}
Q(\bar{a}*\bar{b},E)
 = -\omega(R(a)b+aR(b)+ab,E)
 = -\omega(a,R^*(b)+R(b)+b)=0.
\end{equation}
That is, $Q(\bar{A}*\bar{A},E)=0$. Therefore, $\bar{A}*\bar{A}\neq \bar{A}$ and $\bar{A}$ is not a semisimple algebra. Let $\bar{A}=B+N$, where $N$ is the nil-radical of $\bar{A}$ and $B$ is the semisimple component of~$\bar{A}$. 
Consider $\bar{E}=\bar{E_s}+\bar{E_n}$, where $\bar{E_s}\in B$, $\bar{E_n}\in N$. Note that $\bar{E_s}\in \bar{A}*\bar{A}$ and by \eqref{th6} $Q(E,\bar{E_s})=0$. Thus, $Q(E,\bar{E_n})=Q(E,\bar{E})\neq 0$.

Consider a map $i\colon A\mapsto I(A)$ defined by~\eqref{i-map}. Equality \eqref{u2} implies that $\bar{E_n}*i(E) =i(E)*\bar{E_n} =i(-R(E_n))$. 
Therefore, $\bar{E_n}*i(E)\in I(A)$ and $\bar{E_n}*i(E)$ is nilpotent too.
Hence, $R(E_n)$ is a~nilpotent matrix. 

Let $J$ be the ideal of $D(A)$ defined by~\eqref{J-map} and $j \colon A\to J$ is the map defined as $j(a)=-\bar{a}-R(a)$ for all $a\in A$. Recall that $-E_n=R(E_n)+R^*(E_n)$. Since $i(E_n)+j(E_n)=E_n$, 
\begin{multline*}
j(E)*\bar{E_n}
 = (E-i(E))*\bar{E_n}
 = \bar{E}_n+i(R(E_n))
 = \overline{E_n+R(E_n)}+R(E_n+R(E_n)) \\
 = j(R^*(E_n))
 = \bar{E_n}*j(E).
\end{multline*}

Therefore, $R^*(E_n)$ is a nilpotent element too. Then $\tr(E_n)=-\tr(R(E_n))-\tr(R^*(E_n))=0$, it is a contradiction,
since $\tr(E_n)=\omega(E_n,E)=Q(\bar{E_n},E)\neq 0$. \hfill $\square$

{\bf Remark 5}. 
We have proved that there are no Rota---Baxter operators of weight $\lambda$ on $M_n(F)$ satisfying $R+R^*+\lambda\id=0$. The same result can be proved for quadratic finite-dimensional simple Jordan and alternative algebras, since simple finite dimensional algebras in these varieties are unital. However, everything changes if we consider a simple finite-dimensional Lie (or Malcev) algebra. In \cite{Goncharov2} it was proved that Rota---Baxter operators of nonzero weight~$\lambda$ satisfying $\theta_R=0$ on a simple finite dimensional Lie (or Malcev) algebra $L$ are in one-to-one correspondence with solutions of the modified Yang---Baxter equation on $L$.

\section{$\lambda$-double Lie algebras}

In the light of Theorem~6, we come to the following definition.

{\bf Definition 4}. 
A $\lambda$-double Lie algebra is a linear space $V$ equipped with a double bracket
satisfying the following identities
\begin{gather}
\dbl a,b\dbr
 + \dbl b,a\dbr^{(12)} = \lambda(a\otimes b-b\otimes a),
\label{lambda-antiCom} \\
\dbl a, \dbl b,c\dbr \dbr _L -\dbl b, \dbl a,c\dbr \dbr _R
 - \dbl \dbl a,b\dbr,c\dbr _L = -\lambda (b\otimes \dbl a,c\dbr)^{(12)} \label{lambda-Jacobi}
\end{gather}
for $a,b,c\in V$.

A $\lambda$-double Lie algebra for $\lambda = 0$ is an ordinary double Lie algebra.

Let $(A,\omega)$ be a quadratic algebra with a~unit. Define  $\tr(a):=\omega(a,1)$. 

{\bf Definition 5}. 
Given a quadratic algebra~$A$, a~linear operator~$R$ on~$A$ is called $\lambda$-skew-symmetric if 
\begin{equation}\label{RBLambdaSkewSym}
R(a)+R^*(a)+\lambda a=\lambda \tr(a)1
\end{equation}
for all $a\in A$.

{\bf Remark 6}.
Note that the condition~\eqref{RBLambdaSkewSym} already appeared in the article of O.~Ogievetsky and T.~Popov~\cite{Ogievetsky}.

Let $V$ be a~finite-dimensional vector space with a double bracket $\dbl \cdot,\cdot \dbr$ and let $R$~be the corresponding operator on $\End(V)$, see~\eqref{eq:Bracket_via_RB}. Then the identity~\eqref{lambda-antiCom}
has the form $\theta_R(y) = \lambda\tr(y)E$.
Further, the identity~\eqref{lambda-Jacobi}
is fulfilled modulo~\eqref{lambda-antiCom} if and only if 
$R$ is an RB-operator of weight~$\lambda$ on $\End(V)$. 
So, the following result holds.

{\bf Theorem~7}.
Let $V$~be a finite-dimensional vector space with a double bracket $\dbl \cdot,\cdot \dbr$
determined by an operator $R\colon\End(V)\to \End(V)$ as in~\eqref{eq:Bracket_via_RB}.
Then $V$ is a $\lambda$-double Lie algebra if and only if $R$ is a~$\lambda$-skew-symmetric RB-operator of weight~$\lambda$ on $\End(V)$.

Let us call a solution~$r$ of AYBE of weight~$\lambda$ on $M_n(F)$ as a $\lambda$-skew-symmetric if
$r$~satisfies the equality
$$
r + \tau(r) = \lambda(E\otimes E-C)
$$
where $C = \sum\limits_{i,j=1}^n e_{ij}\otimes e_{ji}$.

Extending~\cite{AYBE-ext}, we get the one-to-one correspondence between $(-\lambda)$-skew-symmetric
solutions of AYBE of weight~$-\lambda$ on $M_n(F)$ and $\lambda$-skew-symmetric RB-operators of weight~$\lambda$ on $M_n(F)$.

{\bf Example 1}.
A linear map $P$ defined on $M_n(F)$ as follows,
$$
P(e_{ij})
 = \begin{cases}
 -e_{ij}, & i<j, \\
 0, & i>j, \\
 \sum\limits_{k\geq1}e_{i+k,\,i+k}, & i=j,
 \end{cases}
$$
is an RB-operator of weight~1 on $M_n(F)$~\cite{Unital}.
By the definition, we have
$$
P^*(e_{ij})
 = \begin{cases}
 0, & i<j, \\
 -e_{ij}, & i>j, \\
 \sum\limits_{k\geq1}e_{i-k,\,i-k}, & i=j,
 \end{cases}
$$
so $P$ is 1-skew-symmetric.
Let $\{f_i\}$ be a basis of $n$-dimensional space~$V$ such that
$e_{ij}(f_k) = \delta_{jk}f_i$.
Due to the equivalence between RB-operators and double brackets, we have
\begin{equation} \label{Ex10}
\dbl f_k,f_l\dbr
 = \begin{cases}
 f_k\otimes f_l - f_l\otimes f_k, & k<l, \\
 0, & k\geq l.
 \end{cases}
\end{equation}

{\bf Remark 7}.
Given a $\lambda$-skew symmetric RB-operator~$R$ of weight~$\lambda$ on $M_n(F)$,
the operator $R^{(T)}$ which acts by the rule $R^{(T)}(a) = (R(a^T))^T$,
where $T$ denotes the transpose in $M_n(F)$,
is again a $\lambda$-skew symmetric RB-operator of weight~$\lambda$ on $M_n(F)$.
Indeed, Proposition~1 implies that $R^{(T)}$ is an RB-operator of weight~$\lambda$.
The operator $R^{(T)}$ is $\lambda$-skew symmetric,
since the property $(S^{(T)})^* = (S^*)^{(T)}$ holds for all operators $S$ on $M_n(F)$.

{\bf Example 2}.
For $P$ from Example~1, consider
$$
P^{(T)}(e_{ij})
 = \begin{cases}
  0, & i<j, \\
  -e_{ij}, & i>j, \\
 \sum\limits_{k\geq1}e_{i+k,\,i+k}, & i=j,
 \end{cases}
$$
it is also a 1-skew-symmetric RB-operator of weight~1 on $M_n(F)$.
Then
\begin{equation} \label{Ex11}
\dbl f_k,f_l\dbr
 =  \begin{cases}
 f_k\otimes f_l, & k<l, \\
 - f_l\otimes f_k, & k>l, \\
 0, & k = l.
 \end{cases}
\end{equation}

The RB-operator from the following example is close to the RB-operator from Example~2 for $n=3$, the only difference is another action on diagonal matrices. So, the obtained double bracket is like a~join of the double brackets~\eqref{Ex10} and~\eqref{Ex11}.

{\bf Example 3}~\cite{Arthamonov}.
Let $A$ be a three-dimensional vector space with a basis $a_1,a_2,a_3$. Define the double bracket on $A\otimes A$:
\begin{equation}\label{Art:Exm}
\begin{gathered}
\dbl a_1,a_2\dbr  = -a_1\otimes a_2,\quad
\dbl a_2,a_1\dbr  =  a_1\otimes a_2,\quad
\dbl a_2,a_3\dbr  =  a_3\otimes a_2,\\
\dbl a_3,a_1\dbr  =  a_1\otimes a_3 - a_3\otimes a_1,\quad
\dbl a_3,a_2\dbr  = -a_3\otimes a_2.
\end{gathered}
\end{equation}
All omitted brackets of generators are assumed to be zero.
It is $(-1)$-double Lie algebra.
In~\cite{Arthamonov}, this double bracket was defined in the context of so called modified Poisson algebra, see the next section. 

The corresponding linear operator $R$ on $M_3(F)$ defined by~\eqref{eq:Bracket_via_RB}
for this double bracket equals
\begin{gather*}
R(e_{12}) = R(e_{13}) = R(e_{23}) = 0, \quad
R(e_{21}) = e_{21}, \quad R(e_{31}) = e_{31}, \quad R(e_{32}) = e_{32}, \\
R(e_{11}) = -e_{22},\quad R(e_{22}) = 0,\quad R(e_{33}) = -(e_{11}+e_{22}),
\end{gather*}
and it is a Rota---Baxter operator of weight~$-1$ on $M_3(F)$ 
(see the case 1) from A)~\cite[Theorem~3]{GonGub}).
Since 
\begin{gather*}
R^*(e_{12}) = e_{12}, \quad 
R^*(e_{13}) = e_{13}, \quad 
R^*(e_{23}) = e_{23}, \quad
R^*(e_{21}) = R^*(e_{31}) = R^*(e_{32}) = 0, \\
R^*(e_{11}) = -e_{33},\quad R^*(e_{22}) = -(e_{11}+e_{33}),\quad R^*(e_{33}) = 0,
\end{gather*}
$R$ is $(-1)$-skew-symmetric.

{\bf Example 4}.
A linear map $P_1$ defined on $M_n(F)$ as follows,
$$
P_1(e_{ij})
 = \begin{cases}
  \sum\limits_{k\geq1}e_{i+k,j+k}, & i\leq j, \\
- \sum\limits_{k\geq0}e_{i-k,j-k}, & i>j,
   \end{cases}
$$
is an RB-operator of weight~1 on $M_n(F)$~\cite{Spectrum}.
Since
$$
P_1^*(e_{ij})
 = \begin{cases}
  \sum\limits_{k\geq1}e_{i-k,j-k}, & i\geq j, \\
- \sum\limits_{k\geq0}e_{i+k,j+k}, & i<j,
  \end{cases}
$$
$P_1$ is 1-skew-symmetric.
Then 
$$
\dbl f_k,f_l\dbr
 = \begin{cases}
 -(f_l\otimes f_k + f_{l+1}\otimes f_{k-1}+\ldots+f_{k-1}\otimes f_{l+1}), & l<k, \\
 f_k\otimes f_l+f_{k+1}\otimes f_{l-1}+\ldots+f_{l-1}\otimes f_{k+1}, & l\geq k,
 \end{cases}
$$
is a 1-double bracket.

If we extend Example~4 on the case of countable-dimensional vector space, we get the space $V = \Bbbk[t]$ equipped with the 1-double Lie bracket
$$
\dbl t^n,t^m\dbr = -\frac{(t^n\otimes t^{m+1}-t^m\otimes t^{n+1})}{t\otimes 1-1\otimes t}.
$$
Here we identify $t^k$ with $f_{k+1}$, $k\geq0$.
Denote the obtained 1-double Lie algebra as $M_1$.

{\bf Remark 8}.
We may prove the analog of Theorem~4 for $\lambda$-double Lie algebras and Rota---Baxter operators from $R\in \End_f(V)$ when $V$ is countable-dimensional.

Let us transform Example~4 as follows.
Define the RB-operator $P_2 = P_1^{(\psi_n)}$ on~$M_n(F)$,
where $\psi_n$ is an automorphism of $M_n(F)$ defined by the formula
$\psi(e_{ij}) = e_{n+1-i,n+1-j}$. Then we extend $P_2$ as an operator from~$I$ to~$\End(V)$.

{\bf Example 5}.
A linear map $P_2\colon I\to \End(V)$ defined as follows,
$$
P_2(e_{ij})
 = \begin{cases}
  \sum\limits_{k\geq1}e_{i-k,j-k}, & i\geq j, \\
- \sum\limits_{k\geq0}e_{i+k,j+k}, & i<j,
   \end{cases}
$$
is a 1-skew-symmetric RB-operator of weight~1 from $I$ to $\End(V)$. 
Then 
$$
\dbl f_k,f_l\dbr
 = \begin{cases}
   f_k\otimes f_l+f_{k-1}\otimes f_{l+1}+\ldots+f_{l+1}\otimes f_{k-1}, & k>l, \\
 -(f_{k+1}\otimes f_{l-1} + f_{k+2}\otimes f_{l-2}+\ldots+f_l\otimes f_k), & k\leq l,
 \end{cases}
$$
is a 1-double bracket.

A vector space $V = \Bbbk[t]$ equipped with~a double bracket
$$
\dbl t^n,t^m\dbr = \frac{t^{n+1}\otimes t^m-t^{m+1}\otimes t^n}{t\otimes 1-1\otimes t}
$$
is a 1-double Lie algebra, denote it as $M_2$.

{\bf Proposition 6}. 
Each of 1-double Lie algebras $M_1$ and $M_2$ 
has only one nonzero proper ideal $I = tF[t]$. Moreover,
$I$ is isomorphic to the whole double Lie algebra.

{\sc Proof}. 
Let us prove the statement for $M_1$, the proof for~$M_2$ is analogous.

Suppose that $I$ is a~nonzero proper ideal in $M_1$.
Define $n$ as the minimal degree of elements from~$I$.
Let us show that $n = 1$ and $t\in I$.

If $n = 0$, then $1\in I$. Let us prove by induction on~$s\geq0$ that $t^s\in I$. For $s = 0$, it is true. Suppose that $s>0$ and we have proved that $t^j\in I$ for all $j<s$. Consider
$$
\dbl 1,t^{2s}\dbr
 = t^{2s-1}\otimes t + t^{2s-2}\otimes t^2 + \ldots + t^{s+1}\otimes t^{s-1}
 + t^s\otimes t^s + \ldots + 1\otimes t^{2s}.
$$
So, $t^s\otimes t^s\in V\otimes I + I\otimes V$.
Consider the map $\psi\colon V\otimes V\to V/I\otimes V/I$
acting as follows, $\psi(v\otimes w) = (v+I)\otimes (w+I)$.
Applying the equality $I\otimes V + V\otimes I = \ker(\psi)$,
we conclude that $\psi(t^s\otimes t^s) = 0$, and it means that $t^s \in I$. 
Thus, $I = M_1$ and it is not proper ideal, as required.

For $n\geq 1$, consider $f = \sum\limits_{j=0}^n\alpha_j t^j\in I$.
We have that the product
$$
\dbl 1,f\dbr
 = \sum\limits_{j=1}^n\alpha_j (t^{j-1}\otimes t + \ldots + t\otimes t^{j-1})
 + 1\otimes f - \alpha_0 1\otimes 1
$$
lies in $V\otimes I + I\otimes V$.
When $n>1$, the elements $1+I,t+I,\ldots,t^{n-1}+I$ of $V/I$ are linearly independent, we obtain a~contradiction. By the same reason, the case $n = 1$ and $\alpha_0\neq0$ does not hold. So, $n = 1$ and $t\in I$.

As above, we may prove by induction on $s\geq1$ that $t^s\in I$. 
For this, it is enough to analyze the double product $\dbl t,t^{2s-1}\dbr$.

Finally, the linear map $\xi\colon M_1\to I$ defined by the formula
$\xi(t^n) = t^{n+1}$ is the isomorphism between $M_1$ and $I$.
\hfill $\square$

{\bf Proposition 7}. 
Let $A$ be a quadratic algebra and let $R$ be a~$\lambda$-skew-symmetric Rota---Baxter operator of weight~$\lambda$ on~$A$. Then $R^*$ is also a~$\lambda$-skew-symmetric Rota---Baxter operator of weight $\lambda$ on $A$.

{\sc Proof}. 
Since $R$ is $\lambda$-skew-symmetric, then 
$R^*(a)=-R(a)-\lambda a+\lambda \tr(a)1$.
Therefore,
\begin{multline}
R^*(a)R^*(b)
 = R(a)R(b)+\lambda R(a) b+\lambda aR(b)-\lambda \tr(b) R(a)-\lambda \tr(a) R(b) \\
+\lambda^2 ab-\lambda^2 \tr(b) a-\lambda^2\tr(a)b+\lambda^2\tr(a)\tr(b)1.
\end{multline}
On the other hand, 
\begin{multline}
R^*(R^*(a)b+aR^*(b)+\lambda ab)=R^*(-R(a)b-\lambda ab+\lambda \tr(a) b-aR(b)-\lambda ab +\lambda \tr(b)a+\lambda ab) \\
= R^*(-R(a)b-aR(b)-\lambda ab)+\lambda R^*(\tr(b)a+\tr(a)b)\\
= R(a)R(b)+\lambda(R(a)b+aR(b)+\lambda ab)-\lambda\tr(R(a)b+aR(b)+\lambda ab)\\
-\lambda R(\tr(b)a+\tr(a)b)-\lambda^2 \tr(a)b-\lambda^2\tr(b)a+2\lambda^2\tr(a)\tr(b).
\end{multline}
Finally,
\begin{multline*}
R^*(a)R^*(b)-R^*(R^*(a)b+aR^*(b)+\lambda ab)
=\lambda\tr(R(a)b+aR(b)+\lambda ab)-\lambda^2\tr(a)\tr(b)\\
=\lambda \tr(R(a)b+aR(b)+\lambda ab-\lambda \tr(a)b)=0,
\end{multline*}
since 
$$
\tr(aR(b))
 = \omega(aR(b),1)
 = \omega(a,R(b))
 = \omega(R^*(a),b)
 = \omega(R^*(a)b,1)
 = \tr(R^*(a)b)
$$
and $R$ is~$\lambda$-skew-symmetric. \hfill $\square$

{\bf Remark 9}. 
Proposition~7 implies that given a quadratic algebra~$(A,\omega)$ equipped with a~$\lambda$-skew-symmetric Rota---Baxter operator of weight~$\lambda$, both linear operators $-R-\lambda\id$ and $-R-\lambda\id+\lambda\omega(\cdot,1)1$ are RB-operators of weight~$\lambda$.

{\bf Remark 10}. 
Let $L$ be a~finite-dimensional $\lambda$-double Lie algebra for nonzero~$\lambda$ and let $R$ be a corresponding $\lambda$-skew-symmetric RB-operator of weight~$\lambda$ on $\End(L)$.
By Proposition~7, we get that $R^*$ is also $\lambda$-skew-symmetric RB-operator on $\End(L)$. Thus, we may define a~new $\lambda$-double Lie algebra structure on the vector space~$L$ by~$R^*$ instead of~$R$.
RB-operators $P_1$ and $P_2$ from Examples~4 and~5 are dual to each other in this sense.

Let $n = \dim(V)>1$, $A = \End(V)$, and $\omega$ is the trace form on $A$.
By~Corollary~1, we have the decomposition
$A = I_1\oplus I_2$ (as subalgebras), where
$$
I_1=\ker(R^{N}),\quad I_2=\ker(R+\lambda\id)^{N}, \quad N=n^2. 
$$
For the Rota---Baxter operator $R^*$ we analogously have the decomposition $A = J_1\oplus J_2$ for $J_1=\ker(R^*)^{N}$ and $J_2=\ker(R^*+\lambda\id)^{N}$. 

Define $I_2' = \ker(R+\lambda\id)$. Let us show that $I_2'\neq(0)$. Suppose that $I_2' = (0)$, then $I_2 = (0)$ and $I_1 = A$. 
Also, $-(R+\lambda\id)$ is an invertible RB-operator of weight~$\lambda$ on $A$. It is well-known that $R+\lambda\id$ is a~homomorphism from $B$ to~$A$~\cite{Splitting}, where $B$ is the vector space $\End(V)$ under the product
$$
x\circ y 
 = -(R+\lambda\id)(x)y - x(R+\lambda\id)(y) + \lambda xy
 = -( R(x)y + xR(y) + \lambda xy ). 
$$
Since $R+\lambda\id$ is nondegenerate, it is an isomorphism between $B$ and $A$.
On the other hand, $\ker(R)$ is a nonzero ideal of $B$ as the kernel of the homomorphism $R\colon B\to A$. Since $B\cong A$~is simple, we conclude that $B = \ker(R)$, i.\,e., $R = 0$. Thus, $R^* = 0$ and $R$~may not be $\lambda$-skew-symmetric.

{\bf Proposition 8}. 
We have

a) $\omega(I_1,J_2)=\omega(J_1,I_2)=0$; 

b) $I_2'$ is a nilpotent ideal in $R(A)$. 

{\sc Proof}. 
a) Let $x\in I_1$, $y\in J_2$. The restriction of the map $R+\lambda\id$ on $I_1$ is nondegenerate. Therefore,  $x=(R+\lambda\id)^N(z)$ for some $z\in I_1$. Then
$$
\omega(x,y)
 = \omega((R+\lambda\id)^N(z),y)
 = \omega(z,(R^*+\lambda\id)^N(y))
 = 0.
$$
Similarly, $\omega(J_1,I_2)=0$.

b) It is well-known that  $I_2'$ is an ideal in $R(A)$, see, e.\,g.,~\cite[Lemma~8]{Spectrum}.

Suppose that $I_2'$ is not nilpotent, then there exists a nonzero idempotent $e^2=e\in I_2'$. Since the trace of $e$ is a positive integer number, $R^*(e)=\lambda\tr(e)E\neq 0$. On the other hand, applying Proposition~7, we get
$$
\lambda^2\tr(e)^2E
 = R^*(e)R^*(e)
 = R^*(2\lambda \tr(e)e + \lambda e)
 = \lambda^2(2\tr(e)^2+\tr(e))E
$$
and $\tr(e)^2+\tr(e)=0$, a~contradiction. \hfill $\square$

{\bf Corollary 2}. 
For every $x\in I_2'$ we have $\tr(x)=0$ and $R^*(x)=0$.

{\bf Corollary 3}. 
There is a natural isomorphism between:

a) $J_1$ and $I_1^*$,

b) $J_2$ and $I_2^*$.

{\sc Proof}.
Indeed, for every $f\in A$ define a map $\gamma: A\mapsto A^*$ as follows:
$$
\gamma(f)(a)=\omega(f,a).
$$
Since the form $\omega$ is non-degenerate, $\gamma$ is an isomorphism. By Proposition 8a, $\gamma(J_1)=I_1^*$ and $\gamma(J_2)=I_2^*$. \hfill $\square$

In particular, if $e_1,\ldots e_p$ is a basis of $I_1$ and $e_{p+1},\ldots, e_N$ is a basis of $I_2$, then we may choose the dual basis $f_1,\ldots f_N$ of $e_1,\ldots e_N$ in such a way that  $f_1,\ldots,f_p\in J_1$ and $f_{p+1},\ldots, f_N\in J_2$.

{\bf Lemma}. 
Let $(V,\dbl\cdot,\cdot\dbr)$ be a $\lambda$-double Lie algebra for $\lambda\neq0$ and let $R$ be the corresponding Rota---Baxter operator of weight~$\lambda$ on $A=\End(V)$. If $\dim V>1$, then $U=I_2'V$ is a~proper ideal of~$V$.

{\sc Proof}.
Let $\lambda=1$. Since $I_2'$ is a nilpotent nonzero ideal in $R(A)$, $U$ is a nonzero proper subspace of $V$. Moreover, $R(A)U\subset U$. It means that for all $a\in V$ and $b\in U$
$$
\dbl a,b\dbr  = \sum\limits_i e_i(a)\otimes R(e_i^*)(b)\in V\otimes U.
$$

It remains to prove the inclusion $\dbl b,a\dbr\in V\otimes U + U\otimes V$ for the same $a,b$. If $y\in J_1$ and $x\in I_2'$, then
\begin{equation} \label{J_1U<U}
yx=-R^*(y)x-R(y)x+\tr(y)x\in -R^*(y)x+I_2'\subset (R^*)^2(y)x+I_2'\subset\ldots\subset I_2'.
\end{equation}
This means that $J_1U\subset U$. 

Let $e_1,\ldots,e_k$ be a basis of $I_2'$, $e_1,\ldots,e_k,e_{k+1},\ldots, e_t$ be the basis of $I_2$. By Corollary~3, we may find a basis $f_1,\ldots,f_t$ of $J_2$ such that $\omega(e_i,f_j)=\delta_{ij}$ for all $i=1,\ldots,t$. 
In particular, for $j=1,\ldots,k$ we have $f_j^*=e_j\in I_2'$, and so
$$
f_j(b)\otimes R(e_j)(a)=-f_j(b)\otimes e_ja\in V\otimes U. 
$$

Let $f_{t+1},\ldots,f_q$ be a~basis of $J_1$.

If $i=1\ldots,k$ and $j=k+1,\ldots,t$, then $\tr(f_je_i)=0$ and moreover, $\tr(f_j I_2') = 0$ by Proposition~8a. Therefore,
\begin{equation} \label{NoSimple:fjei}
R(f_je_i)+R^*(f_je_i)+f_je_i=0.
\end{equation}

Then by Proposition~8a,
$$
\tr(R^*(f_je_i)a)=\tr(f_je_iR(a))=0,
$$
since $e_iR(a)\in I_2'$.

Therefore, $R^*(f_je_i)=0$ and $f_je_i\in I_2'$ by~\eqref{NoSimple:fjei} for all $i=1\ldots,k$ and $j=k+1,\ldots, t$. Consequently, $f_jU\subset U$ and we get
$f_j(b)\otimes R(e_j)(a)\subset U\otimes V$.

Applying~\eqref{J_1U<U}, we finally have that
\begin{multline*}
\dbl b,a \dbr 
= \sum\limits_{j=1}^k f_j(b)\otimes R(f_j^*)(a)
+\sum\limits_{j=k+1}^t f_j(b)\otimes R(f_j^*)(a)
+ \sum\limits_{j=t+1}^q f_j(b) \otimes R(f_j^*)(a) \\
\subset V\otimes U + U\otimes V.
\end{multline*}
Theorem is proved.
\hfill $\square$

{\bf Remark 11}. 
Let $\dim V>1$ and $J_2'=\ker(R^*+\lambda\id)$. Then the subspace $U'=J_2'V$ is also a proper ideal in $(V,\dbl\cdot,\cdot\dbr)$.

{\bf Theorem 8}. 
There are no simple finite-dimensional $\lambda$-double Lie algebras.

{\sc Proof}.
For $\lambda = 0$, it was proved in~\cite{DoubleLie}.
When $\lambda\neq0$ and $\dim V>1$, it follows from Lemma.
Finally, when $\lambda\neq0$ and $\dim V = 1$, it is easy to show that
we have $\dbl V,V\dbr = 0$, so $V$ is not simple too.
\hfill $\square$

\section{Modified double Poisson algebras}

{\bf Definition 6}~\cite{Arthamonov0,Arthamonov}.
Let $A$ be an associative algebra over $F$ with the product $ab=\mu(a\otimes b)$. 
A double bracket $\dbl\cdot,\cdot\dbr$ on $A$ is called 
a modified double Poisson bracket~\cite{Arthamonov0,Arthamonov} 
if the following equalities hold
\begin{gather}
\dbl a,bc \dbr = (b\otimes 1)\dbl a,c\dbr+\dbl a,b\dbr(1\otimes c), \label{a1} \\
\dbl ab,c \dbr = (1\otimes a)\dbl b,c\dbr+\dbl a,c\dbr(b\otimes 1), \label{a2} \\ 
\{a,\{b,c\}\}-\{b,\{a,c\}\}=\{\{a,b\},c\}, \label{a3} \\
\{a,b\}+\{b,a\}=0\ \text{mod}\ [A,A] \label{a4} 
\end{gather}
for all $a,b,c\in A$. Here $\{a,b\}=\mu\circ \dbl a,b\dbr$.

Note that the identities~\eqref{a1} and~\eqref{a2} are exactly the same as both Leibniz rules~\eqref{Leibniz} and~\eqref{LeibnizTwo} fulfilled in double Poisson algebras.
The identities~\eqref{a3} and~\eqref{a4} are weakened versions of anti-commutativity and Jacobi identity.

S. Arthamonov posed two conjectures devoted to modified double Poisson brackets, one of them is the following.

{\bf Conjecture} (S. Arthamonov, 2017~\cite{Arthamonov}).
The bracket defined on $A\otimes A$ by~\eqref{Art:Exm} is a~modified Poisson bracket.

We prove more general result that every $\lambda$-double Lie algebra generates the structure 
of modified double Poisson algebra on the free associative algebra.

{\bf Theorem 9}. 
Let $(V,\dbl \cdot,\cdot\dbr)$ be a finite-dimensional $\lambda$-double Lie algebra with nonzero~$\lambda$. Then equalities~\eqref{a1},~\eqref{a2} define a modified double Poisson algebra on the free associative algebra $A=\As\langle e_1,\ldots,e_n\rangle$, where $e_1,\ldots,e_n$ is a basis of $V$.

{\sc Proof}.
For convenience, we will prove the statement for $\lambda=1$. Take $a=x_1\ldots x_k$, $b=y_1\ldots y_l$, where $x_i,y_j\in \{e_1,\ldots,e_n\}$. Then
$$
\dbl a,b\dbr =\sum_{i=1}^k\sum_{j=1}^l (y_1\ldots y_{j-1}\otimes x_1\ldots x_{i-1})\dbl x_i, y_j\dbr(x_{i+1}\ldots x_{k}\otimes y_{j+1}\ldots y_l).
$$
It is easy to check that $\dbl\cdot,\cdot\dbr$ satisfies~\eqref{a1} and~\eqref{a2}.

For $x,y\in V$ we will use the following notation: 
$$\dbl x, y\dbr=(x,y)_{(1)}\otimes (x,y)_{(2)}.$$
Let us prove  \eqref{a4}. We have
\begin{multline}\label{th**}
\{a,b\}=\sum\limits_{i,j} y_1\ldots y_{j-1}(x_i,y_j)_{(1)}x_{i+1}\ldots x_kx_1\ldots x_{i-1}(x_i,y_j)_{(2)}y_{j+1}\ldots y_l\\
=\sum\limits_{i,j}x_1\ldots x_{i-1}(x_i,y_j)_{(2)}y_{j+1}\ldots y_ly_1\ldots y_{j-1}(x_i,y_j)_{(1)}x_{i+1}\ldots x_k\ \text{mod}\ [A,A].
\end{multline}

The equality \eqref{lambda-antiCom} implies that for all $a,b,c\in A$, 

\begin{equation}\label{tozh}
a(x_i,y_j)_{(2)}b(x_i,y_j)_{(1)}c
 = -a(y_j,x_i)_{(1)}b(y_j,x_i)_{(2)}c+ay_jbx_ic-ax_iby_jc.
\end{equation}
We have
\begin{multline*}
    \sum\limits_{i,j}x_1\ldots x_{i-1}(x_i,y_j)_{(2)}y_{j+1}\ldots y_ly_1\ldots y_{j-1}(x_i,y_j)_{(1)}x_{i+1}\ldots x_k\\
    =-\sum\limits_{i,j}x_1\ldots x_{i-1}(y_j,x_i)_{(1)}y_{j+1}\ldots y_ly_1\ldots y_{j-1}(y_j,x_i)_{(2)}x_{i+1}\ldots x_k \allowdisplaybreaks \\
    +\sum\limits_{i,j}(x_1\ldots x_{i-1}y_jy_{j+1}\ldots y_ly_1\ldots y_{j-1}x_ix_{i+1}\ldots x_k\\
    -x_1\ldots x_{i-1}x_iy_{j+1}\ldots y_ly_1\ldots y_{j-1}y_jx_{i+1}\ldots x_k)\\
    =-\{b,a\}+\sum_{i,j} x_1\ldots x_{i-1}( y_jy_{j+1}\ldots y_ly_1\ldots y_{j-1}x_i-x_iy_{j+1}\ldots y_ly_1\ldots y_{j-1}y_j)x_{i+1}\ldots x_k.
\end{multline*}

Note that
$$
\sum_{j=1}^l y_jy_{j+1}\ldots y_ly_1\ldots y_{j-1}x_i-x_iy_{j+1}\ldots y_ly_1\ldots y_{j-1}y_j
=\left[\sum_j y_{j+1}\ldots y_ly_1\ldots y_{j-1},x_i\right].
$$

Finally,
\begin{multline*}
    \sum\limits_{i=1}^k x_1\ldots x_{i-1}\left(\sum_{j=1}^l y_jy_{j+1}\ldots y_ly_1\ldots y_{j-1}x_i-x_iy_{j+1}\ldots y_ly_1\ldots y_{j-1}y_j\right )x_{i+1}\ldots x_k
    \\
    =\sum_{i=1}^k x_1\ldots x_{i-1}\left[\sum_j y_{j+1}\ldots y_ly_1\ldots y_{j-1},x_i\right]x_{i+1}\ldots x_k\\
    =\sum_{j=1}^l [y_{j+1}\ldots y_ly_1\ldots y_{j-1}, x_1\ldots x_k]\in [A,A].
\end{multline*}

Therefore, $\dbl\cdot,\cdot\dbr$ satisfies \eqref{a4}.

Let us prove~\eqref{a3}.
Define $L(a,b,c) = \{a,\{b,c\}\}-\{b,\{a,c\}\}-\{\{a,b\},c\}$. 
First we note that induction on $\deg(c)$ allows us to assume that $\deg(c) = 1$, 
since for $c = c_1c_2$ we have by~\eqref{a1},~\eqref{a2}:
\begin{multline*}
L(a,b,c)
 = \{a,\{b,c_1c_2\}\}    
 - \{b,\{a,c_1c_2\}\}
 - \{\{a,b\},c_1c_2\} \\
 = \{a,c_1\{b,c_2\}\} + \{a,\{b,c_1\}c_2\}      
 - \{b,c_1\{a,c_2\}\} - \{b,\{a,c_1\}c_2\}
 - \{\{a,b\},c_1c_2\} \\
 = c_1\{a,\{b,c_2\}\} + \underline{\{a,c_1\}\{b,c_2\}}
 + \underline{\underline{\{b,c_1\}\{a,c_2\}}} + \{a,\{b,c_1\}\}c_2 
 - c_1\{b,\{a,c_2\}\} \\
 - \underline{\underline{\{b,c_1\}\{a,c_2\} }}
 - \underline{\{a,c_1\}\{b,c_2\}} - \{b,\{a,c_1\}\}c_2
 - c_1\{\{a,b\},c_2\} - \{\{a,b\},c_1\}c_2 \\
 = c_1L(a,b,c_2) + L(a,b,c_1)c_2. 
\end{multline*}

Let $a=x_1,\ldots,x_k$, $b=y_1,\ldots,y_l$, where $x_i,y_j\in \{e_1,\ldots,e_n\}$ and $c\in V$. We will use the following notation: $\alpha_i(x):=x_1\ldots x_{i-1}$ (for convenience, $\alpha_1(x):=1$), $\beta_i(x)=x_{i+1}\ldots x_k$ ($\beta_k(x):=1)$. That is, $a=\alpha_i(x)x_i\beta_i(x)$. Similarly, $\alpha_i(y)=y_1\ldots y_{i-1}$, $\beta_i(y)=y_{i+1}\ldots y_l$. 
Also, we need
$\gamma_{i,j}(x)
 = \begin{cases}
 x_{i}x_{i+1}\ldots x_{j}, & i\leq j,\\
 1, & i>j.
 \end{cases}$ 

We have
$$
\{b,c\}=\sum\limits_{j=1}^l(y_j,c)_{(1)}\beta_j(y)\alpha_j(y)(y_j,c)_{(2)}.
$$

Therefore,
\begin{align}
    &\{a,\{b,c\}\}  =\sum\limits_{j=1}^l\{a, (y_j,c)_{(1)}\beta_j(y)\alpha_j(y)(y_j,c)_{(2)}\} \nonumber \allowdisplaybreaks \\
    & \ =\sum\limits_{i=1}^k \sum\limits_{j=1}^l(x_i,(y_j,c)_{(1)})_{(1)}\beta_i(x)\alpha_i(x)(x_i,(y_j,c)_{(1)})_{(2)}\beta_j(y)\alpha_j(y)(y_j,c)_{(2)} \nonumber \\
    & \ +\sum\limits_{i=1}^k \sum\limits_{j=1}^{l-1}\sum\limits_{s=j+1}^l(y_j,c)_{(1)}\gamma_{j+1,s-1}(y) (x_i,y_s)_{(1)}\beta_i(x)\alpha_i(x)(x_i,y_s)_{(2)}\beta_s(y)\alpha_j(y)(y_j,c)_{(2)} \label{a(bc)-1} \allowdisplaybreaks \\
    & \ +\sum\limits_{i=1}^k \sum\limits_{j=2}^l\sum\limits_{p=1}^{j-1}(y_j,c)_{(1)}\beta_j(y)\alpha_p(y) (x_i,y_p)_{(1)}\beta_i(x)\alpha_i(x)(x_i,y_p)_{(2)}\gamma_{p+1,j-1}(y)(y_j,c)_{(2)} \label{a(bc)-2} \allowdisplaybreaks \\
& \ +\sum\limits_{i=1}^k \sum\limits_{j=1}^l (y_j,c)_{(1)}\beta_j(y)\alpha_j(y)(x_i,(y_j,c)_{(2)})_{(1)}\beta_i(x)\alpha_i(x)(x_i,(y_j,c)_{(2)})_{(2)}. \nonumber
\end{align}

Similarly,
\begin{align}
    &\{b,\{a,c\}\}=\sum\limits_{i=1}^k\{b, (x_i,c)_{(1)}\beta_i(x)\alpha_i(x)(x_i,c)_{(2)}\} \nonumber \\
    & \ =\sum\limits_{i=1}^k \sum\limits_{j=1}^l(y_j,(x_i,c)_{(1)})_{(1)}\beta_j(y)\alpha_j(y)(y_j,(x_i,c)_{(1)})_{(2)}\beta_i(x)\alpha_i(x)(x_i,c)_{(2)} \nonumber \\
    & \ +\sum\limits_{i=1}^{k-1} \sum\limits_{j=1}^l\sum\limits_{q=i+1}^k(x_i,c)_{(1)}\gamma_{i+1,q-1}(x) (y_j,x_q)_{(1)}\beta_j(y)\alpha_j(y)(y_j,x_q)_{(2)}\beta_q(x)\alpha_i(x)(x_i,c)_{(2)} \label{b(ac)-Psi} \allowdisplaybreaks \\
    & \ +\sum\limits_{i=2}^k \sum\limits_{j=1}^l\sum\limits_{r=1}^{i-1}(x_i,c)_{(1)}\beta_i(x)\alpha_r(x) (y_j,x_r)_{(1)}\beta_j(y)\alpha_j(y)(y_j,x_r)_{(2)}\gamma_{r+1,i-1}(x)(x_i,c)_{(2)} \label{b(ac)-Gamma} \\
& \ +\sum\limits_{i=1}^k \sum\limits_{j=1}^l (x_i,c)_{(1)}\beta_i(x)\alpha_i(x)(y_j,(x_i,c)_{(2)})_{(1)}\beta_j(y)\alpha_j(y)(y_j,(x_i,c)_{(2)})_{(2)}. \nonumber
\end{align}

Finally,
\begin{align}
    & \{\{a,b\},c\}=\sum\limits_{i=1}^k\sum\limits_{j=1}^{l} \{\alpha_j(y)(x_i,y_j)_{(1)}\beta_i(x)\alpha_i(x)(x_i,y_j)_{(2)}\beta_j(y),c\} \nonumber \\
    & \ =\sum\limits_{i=1}^k\sum\limits_{j=2}^{l}\sum\limits_{p=1}^{j-1}(y_p,c)_{(1)}\gamma_{p+1,j-1}(y)(x_i,y_j)_{(1)}\beta_i(x)\alpha_i(x)(x_i,y_j)_{(2)}\beta_j(y)\alpha_p(y)(y_p,c)_{(2)} \label{(ab)c)-1} \\
    & \ +\sum\limits_{i=1}^{k}\sum\limits_{j=1}^{l}((x_i,y_j)_{(1)},c)_{(1)}\beta_i(x)\alpha_i(x)(x_i,y_j)_{(2)}\beta_j(y)\alpha_j(y)((x_i,y_j)_{(1)}),c)_{(2)} \nonumber \allowdisplaybreaks \\
    & \ +\sum\limits_{i=1}^{k-1}\sum\limits_{j=1}^{l}\sum\limits_{q=i+1}^k(x_q,c)_{(1)}\beta_q(x)\alpha_i(x)(x_i,y_j)_{(2)}\beta_j(y)\alpha_j(y)(x_i,y_j)_{(1)}\gamma_{i+1,q-1}(x)(x_q,c)_{(2)} \label{(ab)c-Gamma}  \\
    & \ +\sum\limits_{i=2}^{k}\sum\limits_{j=1}^{l}\sum\limits_{r=1}^{i-1}(x_r,c)_{(1)}\gamma_{r+1,i-1}(x)(x_i,y_j)_{(2)}\beta_j(y)\alpha_j(y)(x_i,y_j)_{(1)}\beta_i(x)\alpha_r(x)(x_r,c)_{(2)} \label{(ab)c-Psi}  \allowdisplaybreaks  \\
    & \ +\sum\limits_{i=1}^k\sum\limits_{j=1}^{l} ((x_i,y_j)_{(2)},c)_{(1)}\beta_j(y)\alpha_j(y)(x_i,y_j)_{(1)}\beta_i(x)\alpha_i(x)((x_i,y_j)_{(2)},c)_{(2)} \nonumber \\
    & \ +\sum\limits_{i=1}^k\sum\limits_{j=1}^{l-1}\sum\limits_{s=j+1}^l\! (y_s,c)_{(1)}\beta_s(y)\alpha_j(y)(x_i,y_j)_{(1)}\beta_i(x)\alpha_i(x)(x_i,y_j)_{(2)}\gamma_{j+1,s-1}(y)(y_s,c)_{(2)} \label{(ab)c-2}. 
\end{align}

First we note that~\eqref{a(bc)-1} with~\eqref{(ab)c)-1}
and~\eqref{a(bc)-2} with~\eqref{(ab)c-2} respectively annihilate.  


By the same arguments and by~\eqref{tozh}, we rewrite the sum of~\eqref{b(ac)-Psi} and~\eqref{(ab)c-Psi} as follows,
\begin{multline*}
    \sum\limits_{i=1}^{k-1} \sum\limits_{j=1}^l\sum\limits_{q=i+1}^k(x_i,c)_{(1)}\gamma_{i+1,q-1}(x) (y_j,x_q)_{(1)}\beta_j(y)\alpha_j(y)(y_j,x_q)_{(2)}\beta_q(x)\alpha_i(x)(x_i,c)_{(2)}\\
    +\sum\limits_{i=2}^{k}\sum\limits_{j=1}^{l}\sum\limits_{r=1}^{i-1}(x_r,c)_{(1)}\gamma_{r+1,i-1}(x)(x_i,y_j)_{(2)}\beta_j(y)\alpha_j(y)(x_i,y_j)_{(1)}\beta_i(x)\alpha_r(x)(x_r,c)_{(2)}
    \\=\sum\limits_{i=1}^{k-1} \sum\limits_{j=1}^l\sum\limits_{q=i+1}^k\left ((x_i,c)_{(1)}\gamma_{i+1,q-1}(x) (y_j,x_q)_{(1)}\beta_j(y)\alpha_j(y)(y_j,x_q)_{(2)}\beta_q(x)\alpha_i(x)(x_i,c)_{(2)}\right.\\
    +\left.(x_i,c)_{(1)}\gamma_{i+1,q-1}(x)(x_q,y_j)_{(2)}\beta_j(y)\alpha_j(y)(x_q,y_j)_{(1)}\beta_q(x)\alpha_i(x)(x_i,c)_{(2)} \right)\\
    =\sum\limits_{i=1}^{k-1} \sum\limits_{j=1}^l\sum\limits_{q=i+1}^k\big((x_i,c)_{(1)}\gamma_{i+1,q-1}(x) y_j\beta_j(y)\alpha_j(y)x_q\beta_q(x)\alpha_i(x)(x_i,c)_{(2)} \allowdisplaybreaks \\
    -\sum\limits_{i=1}^{k-1} \sum\limits_{j=1}^l\sum\limits_{q=i+1}^k ((x_i,c)_{(1)}\gamma_{i+1,q-1}(x) x_q\beta_j(y)\alpha_j(y)y_j\beta_q(x)\alpha_i(x)(x_i,c)_{(2)}\big)=:\Psi.
\end{multline*}

Note that $\sum\limits_{j=1}^ly_j\beta_j(y)\alpha_j(y)=\sum\limits_{j=1}^l\beta_j(y)\alpha_j(y)y_j$. Let $\psi(y)=\sum\limits_{j=1}^l\beta_j(y)\alpha_j(y)y_j$. Then the last equality has a form
\begin{multline*}
    \Psi=\sum\limits_{i=1}^{k-1} \sum\limits_{q=i+1}^k(x_i,c)_{(1)}\gamma_{i+1,q-1}(x) \psi(y)x_q\beta_q(x)\alpha_i(x)(x_i,c)_{(2)}\\
    -\sum\limits_{i=1}^{k-1} \sum\limits_{q=i+1}^k (x_ic)_{(1)}\gamma_{i+1,q-1}(x) x_q\psi(y)\beta_q(x)\alpha_i(x)(x_i,c)_{(2)} \allowdisplaybreaks \\
    =\sum\limits_{i=1}^{k-1}(x_i,c)_{(1)}\psi(y)\beta_i(x)\alpha_i(x)(x_i,c)_{(2)}-(x_i,c)_{(1)}\beta_i(x)\psi(y)\alpha_i(x)(x_i,c)_{(2)}.
\end{multline*}

Similarly, we rewrite the sum of~\eqref{b(ac)-Gamma} and~\eqref{(ab)c-Gamma}
\begin{multline*}
    \sum\limits_{i=2}^k \sum\limits_{j=1}^l\sum\limits_{r=1}^{i-1}(x_i,c)_{(1)}\beta_i(x)\alpha_r(x) (y_j,x_r)_{(1)}\beta_j(y)\alpha_j(y)(y_j,x_r)_{(2)}\gamma_{r+1,i-1}(x)(x_i,c)_{(2)}\\
    +\sum\limits_{i=1}^{k-1}\sum\limits_{j=1}^{l}\sum\limits_{q=i+1}^k(x_q,c)_{(1)}\beta_q(x)\alpha_i(x)(x_i,y_j)_{(2)}\beta_j(y)\alpha_j(y)(x_i,y_j)_{(1)}\gamma_{i+1,q-1}(x)(x_q,c)_{(2)}  \allowdisplaybreaks \\
    =\sum\limits_{i=2}^{k}(x_i,c)_{(1)}\beta_i(x)\psi(y)\alpha_i(x)(x_i,c)_{(2)}-(x_i,c)_{(1)}\beta_i(x)\alpha_i(x)\psi(y)(x_i,c)_{(2)} =: \Gamma.
\end{multline*}

Observe that
\begin{multline*}
\Psi+\Gamma
 = \sum\limits_{i=1}^{k-1}(x_i,c)_{(1)}\psi(y)\beta_i(x)\alpha_i(x)(x_i,c)_{(2)}-(x_i,c)_{(1)}\beta_i(x)\psi(y)\alpha_i(x)(x_i,c)_{(2)}\\
 +\sum\limits_{i=2}^{k}(x_i,c)_{(1)}\beta_i(x)\psi(y)\alpha_i(x)(x_i,c)_{(2)}-(x_i,c)_{(1)}\beta_i(x)\alpha_i(x)\psi(y)(x_i,c)_{(2)} \allowdisplaybreaks \\
 =\sum\limits_{i=1}^{k-1}(x_i,c)_{(1)}\psi(y)\beta_i(x)\alpha_i(x)(x_i,c)_{(2)}-(x_1,c)_{(1)}\beta_1(x)\psi(y)\alpha_1(x)(x_1,c)_{(2)}
 \\
 -\sum\limits_{i=2}^{k-1}(x_i,c)_{(1)}\beta_i(x)\psi(y)\alpha_i(x)(x_i,c)_{(2)}+\sum\limits_{i=2}^{k-1}(x_i,c)_{(1)}\beta_i(x)\psi(y)\alpha_i(x)(x_i,c)_{(2)} \\
 +(x_k,c)_{(1)}\beta_k(x)\psi(y)\alpha_k(x)(x_k,c)_{(2)}-\sum\limits_{i=2}^{k}(x_i,c)_{(1)}\beta_i(x)\alpha_i(x)\psi(y)(x_i,c)_{(2)}\\
=\sum\limits_{i=1}^k (x_i,c)\psi(y)\beta_i(x)\alpha_i(x)(x_i,c)_{(2)}
 - \sum\limits_{i=1}^k(x_i,c)_{(1)}\beta_i(x)\alpha_i(x)\psi(y)(x_i,c)_{(2)}. 
\end{multline*}
 
Summing up the obtained equations we get that
\begin{multline*}
    \{a,\{b,c\}\}-\{b,\{a,c\}\}-\{\{a,b\},c\}\\
    =\sum\limits_{i=1}^k \sum\limits_{j=1}^l(x_i,(y_j,c)_{(1)})_{(1)}\beta_i(x)\alpha_i(x)(x_i,(y_j,c)_{(1)})_{(2)}\beta_j(y)\alpha_j(y)(y_j,c)_{(2)}\\
+\sum\limits_{i=1}^k \sum\limits_{j=1}^l (y_j,c)_{(1)}\beta_j(y)\alpha_j(y)(x_i,(y_j,c)_{(2)})_{(1)}\beta_i(x)\alpha_i(x)(x_i,(y_j,c)_{(2)})_{(2)}\\
-\sum\limits_{i=1}^k \sum\limits_{j=1}^l(y_j,(x_i,c)_{(1)})_{(1)}\beta_j(y)\alpha_j(y)(y_j,(x_i,c)_{(1)})_{(2)}\beta_i(x)\alpha_i(x)(x_i,c)_{(2)}
\\
-\sum\limits_{i=1}^k \sum\limits_{j=1}^l (x_i,c)_{(1)}\beta_i(x)\alpha_i(x)(y_j,(x_i,c)_{(2)})_{(1)}\beta_j(y)\alpha_j(y)(y_j,(x_i,c)_{(2)})_{(2)}  \allowdisplaybreaks \\
-\sum\limits_{i=1}^{k}\sum\limits_{j=1}^{l}((x_i,y_j)_{(1)},c)_{(1)}\beta_i(x)\alpha_i(x)(x_i,y_j)_{(2)}\beta_j(y)\alpha_j(y)((x_i,y_j)_{(1)}),c)_{(2)} \\
-\sum\limits_{i=1}^k\sum\limits_{j=1}^{l} ((x_i,y_j)_{(2)}c)_{(1)}\beta_j(y)\alpha_j(y)(x_i,y_j)_{(1)}\beta_i(x)\alpha_i(x)((x_i,y_j)_{(2)},c)_{(2)}  \\
-\sum\limits_{i=1}^k (x_i,c)\psi(y)\beta_i(x)\alpha_i(x)(x_i,c)_{(2)}+\sum\limits_{i=1}^k(x_i,c)_{(1)}\beta_i(x)\alpha_i(x)\psi(y)(x_i,c)_{(2)}.
\end{multline*}

Note that the sum $-\sum\limits_{i=1}^k\sum\limits_{j=1}^{l} ((x_i,y_j)_{(2)}c)_{(1)}\beta_j(y)\alpha_j(y)(x_i,y_j)_{(1)}\beta_i(x)\alpha_i(x)((x_i,y_j)_{(2)},c)_{(2)}$ can be considered as a linear function on $\dbl x_i,y_j\dbr^{(12)}=(x_i,y_j)_{(2)}\otimes (x_i,y_j)_{(1)}$:
\begin{multline*}
    -\sum\limits_{i=1}^k\sum\limits_{j=1}^{l} ((x_i,y_j)_{(2)}c)_{(1)}\beta_j(y)\alpha_j(y)(x_i,y_j)_{(1)}\beta_i(x)\alpha_i(x)((x_i,y_j)_{(2)},c)_{(2)}\\
    =-\mu^2\left ((1\otimes \beta_j(y)\alpha_j(y)\otimes \beta_i(x)\alpha_i(x))\Theta_c( (x_i,y_j)_{(2)}\otimes (x_i,y_j)_{(1)})\right ),
\end{multline*}
where the map $\mu^2:A\otimes A\otimes A\mapsto A$ is defined as $\mu^2(x\otimes y\otimes z)=xyz$ and
$$
\Theta_c(x\otimes y)=(x,c)_{(1)}\otimes y\otimes (x,c)_{(2)}.
$$

Therefore, we apply \eqref{tozh} and get
\begin{multline*}
    -\sum\limits_{i=1}^k\sum\limits_{j=1}^{l} ((x_i,y_j)_{(2)},c)_{(1)}\beta_j(y)\alpha_j(y)(x_i,y_j)_{(1)}\beta_i(x)\alpha_i(x)((x_i,y_j)_{(2)},c)_{(2)}  \\
    =\sum\limits_{i=1}^k\sum\limits_{j=1}^{l} ((y_j,x_i)_{(1)},c)_{(1)}\beta_j(y)\alpha_j(y)(y_j,x_i)_{(2)}\beta_i(x)\alpha_i(x)((y_j,x_i)_{(1)},c)_{(2)} \\
    -\sum\limits_{i=1}^k\sum\limits_{j=1}^l(y_j,c)_{(1)}\beta_j(y)\alpha_j(y)x_i\beta_i(x)\alpha_i(x)(y_j,c)_{(2)} \allowdisplaybreaks \\
    +\sum\limits_{i=1}^k\sum\limits_{j=1}^l(x_i,c)_{(1)}\beta_j(y)\alpha_j(y)y_j\beta_i(x)\alpha_i(x)(x_i,c)_{(2)}.
\end{multline*}

Let us divide summands from $ \{a,\{b,c\}\}-\{b,\{a,c\}\}-\{\{a,b\},c\}$  into two groups:
$$
\{a,\{b,c\}\}-\{b,\{a,c\}\}-\{\{a,b\},c\} = I+J,
$$
where
\begin{multline*}
I=\sum\limits_{i=1}^k \sum\limits_{j=1}^l(x_i(y_j,c)_{(1)})_{(1)}\beta_i(x)\alpha_i(x)(x_i,(y_j,c)_{(1)})_{(2)}\beta_j(y)\alpha_j(y)(y_j,c)_{(2)}\allowdisplaybreaks \\
-\sum\limits_{i=1}^k \sum\limits_{j=1}^l (x_i,c)_{(1)}\beta_i(x)\alpha_i(x)(y_j,(x_i,c)_{(2)})_{(1)}\beta_j(y)\alpha_j(y)(y_j,(x_i,c)_{(2)})_{(2)} \\
-\sum\limits_{i=1}^{k}\sum\limits_{j=1}^{l}((x_i,y_j)_{(1)},c)_{(1)}\beta_i(x)\alpha_i(x)(x_i,y_j)_{(2)}\beta_j(y)\alpha_j(y)((x_i,y_j)_{(1)},c)_{(2)}\\
+\sum\limits_{i=1}^k(x_i,c)_{(1)}\beta_i(x)\alpha_i(x)\psi(y)(x_i,c)_{(2)},
\end{multline*}

\vspace{-0.5cm}
\begin{multline*} \allowdisplaybreaks
J =\sum\limits_{i=1}^k \sum\limits_{j=1}^l (y_j,c)_{(1)}\beta_j(y)\alpha_j(y)(x_i,(y_j,c)_{(2)})_{(1)}\beta_i(x)\alpha_i(x)(x_i,(y_j,c)_{(2)})_{(2)} \\
    -\sum\limits_{i=1}^k \sum\limits_{j=1}^l(y_j,(x_i,c)_{(1)})_{(1)}\beta_j(y)\alpha_j(y)(y_j,(x_i,c)_{(1)})_{(2)}\beta_i(x)\alpha_i(x)(x_i,c)_{(2)}\\
+\sum\limits_{i=1}^k\sum\limits_{j=1}^{l} ((y_j,x_i)_{(1)},c)_{(1)}\beta_j(y)\alpha_j(y)(y_j,x_i)_{(2)}\beta_i(x)\alpha_i(x)((y_j,x_i)_{(1)},c)_{(2)}\\
-\sum\limits_{i=1}^k\sum\limits_{j=1}^l(y_j,c)_{(1)}\beta_j(y)\alpha_j(y)x_i\beta_i(x)\alpha_i(x)(y_j,c)_{(2)}.
\end{multline*}

Finally, we can use \eqref{lambda-Jacobi} and get that $I=J=0$. Indeed, from \eqref{lambda-Jacobi} it follows that
$$
\dbl x_i,\dbl y_j,c\dbr \dbr_L-\dbl y_j,\dbl x_i,c\dbr\dbr_R-\dbl\dbl x_i,y_j\dbr,c\dbr_L =-(y_j\otimes\dbl x_i,c\dbr)^{(12)}
$$
can be rewritten as
\begin{multline*}
 T_{ij}:= (x_i,(y_j,c)_{(1)})_{(1)}\otimes (x_i,(y_j,c)_{(1)})_{(2)}\otimes (y_j,c)_{(2)} \\
 - (x_i,c)_{(1)}\otimes (y_j,(x_i,c)_{(2)})_{(1)} \otimes (y_j,(x_i,c)_{(2)})_{(2)} \\
 - ((x_i,y_j)_{(1)},c)_{(1)}\otimes (x_i,y_j)_{(2)}\otimes ((x_i,y_j)_{(1)},c)_{(2)} \\ 
 + (x_i,c)_{(1)}\otimes y_j\otimes (x_i,c)_{(2)}
 = 0.
\end{multline*}
Therefore, we get that
$$
\sum\limits_{i=1}^k\sum\limits_{j=1}^lT_{ij}\cdot(\beta_i(x)\alpha_i(x)\otimes \beta_j(y)\alpha_j(y)\otimes 1)=0.
$$
It is left to note that 
$$I=\mu^2\left (\sum\limits_{i=1}^k\sum\limits_{j=1}^lT_{ij}\cdot(\beta_i(x)\alpha_i(x)\otimes \beta_j(y)\alpha_j(y)\otimes 1)\right )=0.$$
The  proof of the equality $J=0$ is similar.
\hfill $\square$

{\bf Corollary 4}.
The double bracket from Example~3 defines a modified double Poisson structure on the algebra $\As\langle a_1,a_2,a_3\rangle$.

Thus, the conjecture of S. Arthamonov~\cite{Arthamonov} holds.

Theorem 9 combined with Corollary 35 from~\cite{Arthamonov} implies the following result.

{\bf Theorem 10}. 
Let $(V,\dbl\cdot,\cdot\dbr)$ be a~finite-dimensional double Lie algebra of nonzero weight~$\lambda$. 
Denote by $\dbl\cdot,\cdot\dbr$ its extension as a modified double Poisson bracket on $\As(V)$. 
Then $F[\Rep_n(\As(V))]^\tr$ is equipped with a Poisson bracket uniquely defined by 
$\{\tr(a), \tr(b)\} = \tr(\dbl a,b\dbr_{(1)}\dbl a,b\dbr_{(2)})$ for any $a, b \in \As(V)$.

\section*{Acknowledgments}

The authors are grateful to the anonymous reviewer for the helpful remarks.

The research is supported by Russian Science Foundation (project 21-11-00286).

\bigskip

\noindent Maxim Goncharov \\
Sobolev Institute of Mathematics \\
Acad. Koptyug ave. 4, 630090 Novosibirsk, Russia \\
Novosibirsk State University \\
Pirogova str. 2, 630090 Novosibirsk, Russia \\
e-mail: gme@math.nsc.ru 

\medskip

\noindent Vsevolod Gubarev \\
Sobolev Institute of Mathematics \\
Novosibirsk State University \\
e-mail: wsewolod89@gmail.com
\end{document}